\documentclass[11pt,a4paper]{article}

\usepackage{graphics,epsfig,theorem,latexsym,amssymb,amsmath, amsfonts}
\usepackage{times,epic,eepic}

\usepackage{pifont}

\usepackage{cite}
\usepackage{url}

\usepackage{color}

\usepackage{newlfont}

 \newenvironment{keywords}{ \noindent {\small\bf Key Words}:}{ }

\def\G1{\hbox{$\displaystyle{\mbox{\ding{172}}}$}}

\def\bd{\begin{description}}
\def\ed{\end{description}}

\def\beq{\begin{equation}}
\def\eeq{\end{equation}}
\def\bea{\begin{eqnarray}}
\def\eea{\end{eqnarray}}
\def\beas{\begin{eqnarray*}}
\def\eeas{\end{eqnarray*}}

\theoremstyle{remark}

\begin{document}

\title{ The exact (up to infinitesimals) infinite perimeter\\ of
the Koch snowflake and its finite area }


\author{    \bf Yaroslav D. Sergeyev\thanks{The author thanks the anonymous reviewers for their
 useful comments. The preparation of the revised version of this paper   was supported
by the Russian Science Foundation, project No.15-11-30022 ``Global
optimization, supercomputing computations, and applications".
}\,\,\,\thanks{Distinguished Professor,
 University of Calabria, Italy; he works also at the
Lobachevsky State University of Nizhni Novgorod,  Russia and the
Institute of High Performance
  Computing and Networking of the National Research Council of
  Italy,         {\tt yaro@si.dimes.unical.it http://wwwinfo.dimes.unical.it/$\sim$yaro   }}
}


\date{}

\maketitle

\vspace*{-5mm}
 \begin{abstract}
The   Koch snowflake   is one of the first fractals that were
mathematically described. It  is  interesting because it has an
infinite perimeter in the limit but its limit area is finite. In
this paper, a recently proposed computational methodology allowing
one to execute numerical computations with infinities and
infinitesimals is applied to study the Koch snowflake at infinity.
Numerical computations with actual infinite and infinitesimal
numbers can be executed on the Infinity Computer being a new
supercomputer patented in USA and EU. It is revealed in the paper
that at infinity the snowflake is not unique, i.e., different
snowflakes can be distinguished   for different infinite numbers of
steps executed during the process of their generation. It is then
shown that for any given infinite number~$n$ of steps it becomes
possible to calculate the exact infinite number, $N_n$, of sides of
the snowflake, the exact infinitesimal length, $L_n$, of each side
and the exact infinite perimeter, $P_n$, of the Koch snowflake as
the result of multiplication of the infinite $N_n$ by the
infinitesimal $L_n$. It is established that for different infinite
$n$ and $k$ the infinite perimeters $P_n$ and $P_k$ are also
different and the difference can be infinite. It is shown that the
finite areas $A_n$ and $A_k$ of the snowflakes can be also
calculated exactly (up to infinitesimals) for different infinite $n$
and $k$ and   the difference $A_n - A_k$ results to be
infinitesimal. Finally, snowflakes constructed starting from
different initial conditions are also studied and their quantitative
characteristics at infinity are computed.

 \end{abstract}


\keywords{ Koch snowflake, fractals, infinite perimeter, finite
area, numerical infinities and infinitesimals, supercomputing.}



\section{Introduction}
\label{s_m1}

Nowadays many fractals are known and their presence  can be found in
nature, especially in physics and biology, in science, and in
engineering (see, e.g.,
\cite{Falconer,Hastings_Sugihara,fractals,fractals2} and references
given therein). Even though   fractal structures were ever around us
their active study started rather recently. The first fractal curves
have been proposed  at the end of XIX$^{th}$ century (see, e.g.,
historical reviews in \cite{Sagan,Sergeyev_Lera_book}) and the word
\textit{fractal} has been introduced by Mandelbrot (see
\cite{Mandelbrot_1,Mandelbrot_2}) in the  second half of the
XX$^{th}$ century. The main geometric characterization of simple
fractals is their self-similarity repeated infinitely many times:
fractals are made by an infinite generation of an increasing number
of smaller and smaller copies of a basic figure  often called an
initiator. More generally, fractal objects  need not exhibit exactly
the same structure at all scales, variations of initiators and
generating procedures (see, e.g., L-systems in \cite{L_systems}) are
conceded. Another important feature of fractal objects is that they
often exhibit fractional dimensions.

\begin{figure}[t]
  \begin{center}
    \epsfig{ figure = 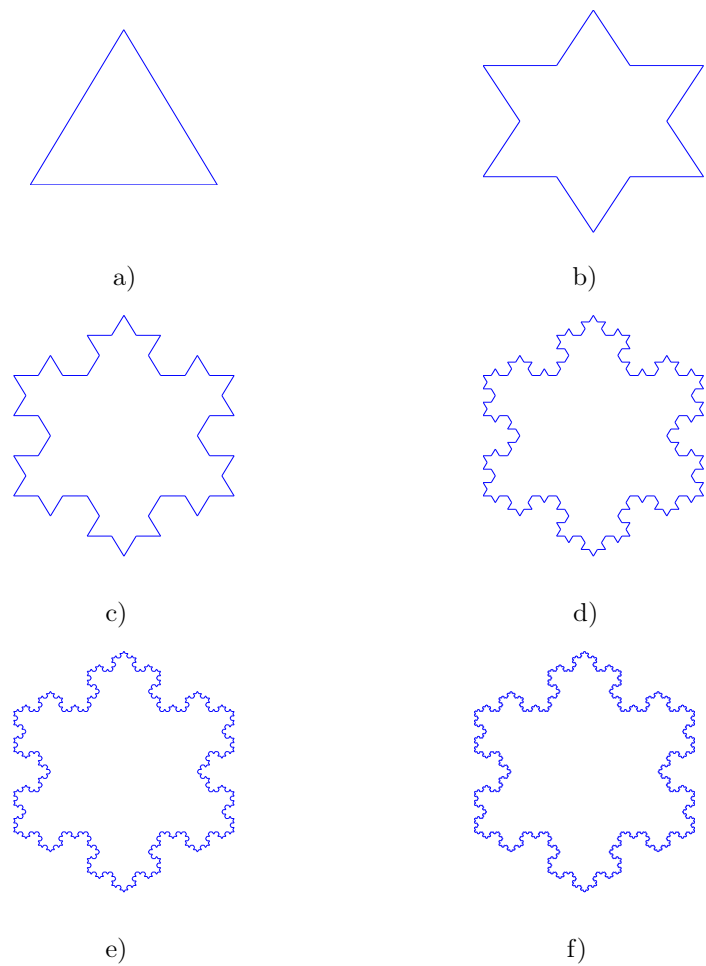, width = 3.6in, height = 5in,  silent = yes }
    \caption{Generation of the Koch snowflake.}
 \label{Koch_fig}
  \end{center}
\end{figure}

Since fractals are objects defined as a limit  of an infinite
process, the computation of their dimension is one of a very few
quantitative characteristics that can be calculated   at infinity.
After $n$ iterations of a fractal process we can give numerical
answers to questions regarding fractals (calculation of, e.g., their
length, area, volume or the number of smaller copies of initiators
present at the $n$-th iteration) only for \textit{finite} values of
$n$. The same questions very often remain without any answer when we
consider an infinite number of steps because when we speak about
limit fractal objects the required values often either tend to zero
and disappear, in practice, or tend to infinity, i.e., become
intractable numerically. Moreover, we cannot distinguish at infinity
fractals starting from similar initiators even though they are
different for any fixed finite value of the iteration number $n$.

In this paper, we study the Koch snowflake that is one of the first
mathematically described fractals. It has been introduced by Helge
von Koch in 1904 (see \cite{Koch}). This fractal is interesting
because it is known that in the limit it has an infinite perimeter
 but its area is finite.
The procedure of its construction is shown in Fig.~\ref{Koch_fig}.
The initiator (iteration number $n=0$) is the triangle shown in
Fig.~\ref{Koch_fig}, a), more precisely, its three sides. Then each
side (segment) is substituted by four smaller segments as it is
shown in Fig.~\ref{Koch_fig}, b) during the first iteration, in
other words, a smaller copy of the triangle is added to each side.
At the iteration $n=2$ (see Fig.~\ref{Koch_fig}, c)) each segment is
substituted again by four smaller segments, an so on. Thus, the Koch
snowflake is the resulting limit object obtained at $n \rightarrow
\infty$. It is known (see, e.g., \cite{fractals}) that its fractal
dimension is equal to $\frac{\log 4}{\log 3}\approx 1.26186$.

Clearly, for finite values of $n$ we can calculate the perimeter
$P_n$ and the respective area $A_n$ of the snowflake. If iteration
numbers $n$ and $k$ are such that $n \neq k$ then it follows  $P_n
\neq  P_k$ and  $A_n \neq  A_k$. Moreover, the snowflake started
from the initiator a) after $n$ iterations will be different with
respect to the snowflake started from the configuration b) after $n$
iterations. The simple illustration for $n=3$ can be viewed in
Fig.~\ref{Koch_fig}. In fact, starting from the initiator a) after
three iterations we have the snowflake d) and starting from the
initiator b) after three iterations we have the snowflake e).

Unfortunately, the traditional analysis of fractals does not allow
us to have quantitative answers to the questions stated above when
$n \rightarrow \infty$. In fact, we know only (see, e.g.,
\cite{fractals}) that
 \beq
\lim_{n \rightarrow \infty} P_n = \lim_{n \rightarrow \infty}
\frac{4^n}{3^{n-1}}l = \infty, \hspace{1cm} \lim_{n \rightarrow
\infty} A_n = \frac{8}{5}a_0, \hspace{1cm} a_0=
\frac{\sqrt{3}}{4}l^2, \label{Koch_1}
 \eeq
where $a_0$ is the area of the original triangle from
Fig.~\ref{Koch_fig}, a) expressed in the terms of its side length
$l$.

In this paper, by using a new computational methodology introduced
in \cite{Sergeyev,informatica,www} we show that    a more precise
quantitative analysis of the Koch snowflake can be done at infinity.
In particular, it becomes possible:
 \bd
 \item[\hspace{6mm}--]
to show that at infinity the Koch snowflake is not a unique object,
namely, different snowflakes can be distinguished at infinity
  similarly to
  different snowflakes that can be distinguished for
different finite values of $n$;
 \item[\hspace{6mm}--]
to calculate the exact (up to infinitesimals) perimeter $P_n$ of the
snowflake (together with the infinite number of sides and the
infinitesimal length of each side) after $n$ iterations for
different \textit{infinite} values of  $n$;
 \item[\hspace{6mm}--]
 to show that for infinite $n$ and $k$ such that $k > n$ it follows
that both $P_n$ and $P_k$ are infinite but $P_k >  P_n$ and their
difference $P_k -  P_n$ can be computed exactly and it results to be
infinite;
\item[\hspace{6mm}--]
to calculate the exact (up to infinitesimals) finite areas $A_n$ and
$A_k$ of the snowflakes   for    infinite $n$ and $k$ such that $k >
n$ and to show that it follows $A_k >  A_n$, the  difference $A_k -
A_n$ can be also calculated exactly and it results to be
infinitesimal;
  \item[\hspace{6mm}--]
to show that the snowflakes constructed starting from different
initiators (e.g., from  initiators shown in Fig.~\ref{Koch_fig}) are
different after $k$ iterations where $k$ is an infinite number,
i.e., they have different infinite perimeters and different areas
where the difference can be measured using infinities and
infinitesimals. \ed

\section{A new numeral system expressing infinities and\\ infinitesimals with a high accuracy}

 In order to understand how it
is possible to study fractals at infinity with an  accuracy that is
higher than that one provided by expressions in (\ref{Koch_1}), let
us remind the important difference that exists between numbers and
numerals. A \textit{numeral} is a symbol or a group of symbols used
to represent  a \textit{number}. A number  is a concept that a
numeral expresses and the difference between them is the same as the
difference between words written in a language and the things the
words refer to. Obviously, the same number can be represented in a
variety of ways by different numerals. For example, the symbols
`11', `eleven', `IIIIIIIIIII',`XI',   and `$\doteq$'   are different
numerals\footnote{The last numeral, $\doteq$, is probably less
known. It belongs to the Maya numeral system where one horizontal
line indicates five and two lines one above the other indicate ten.
Dots are added above the lines to represent additional units. So,
$\doteq$ means eleven and is written as 5+5+1.}, but they all
represent the same number.   A \textit{numeral system} consists of a
set of rules used for writing down numerals and algorithms for
executing arithmetical operations with these numerals. It should be
stressed  that the algorithms  can  vary significantly   in
different numeral systems and their complexity can be also
dissimilar. For instance, division in Roman numerals is extremely
laborious and in the positional numeral system it is much easier.

Notice also that    different  numeral systems  can express
different sets of numbers. One of the simplest existing numeral
systems  that allows its users to express very few numbers is the
system used by Warlpiri people, aborigines living in the Northern
Territory of Australia (see \cite{Butterworth}) and by Pirah\~{a}
people living in Amazonia (see \cite{Gordon}). Both peoples   use
the same very poor numeral system for counting consisting just of
three numerals -- one, two, and `many' -- where `many' is used for
all quantities larger than two.

As a result,   this poor   numeral system does not allow Warlpiri
and Pirah\~{a}   to distinguish numbers larger than 2, to execute
arithmetical operations with them, and, in general, to say a word
about these quantities because in their languages there are neither
words nor concepts for them. In particular, results of operations
2+1 and 2+2 are not   3 and 4 but just `many' since they do not know
about the existence of 3 and 4. It is worthy to emphasize thereupon
that the result `many' is not wrong, it is correct but \textit{its
accuracy is low}. Analogously, when we look at a mob, then both
phrases `There are 2053 persons in the mob' and `There are many
persons in the mob' are correct but the accuracy of the former
phrase  is higher than the accuracy of the latter one.

  Our interest to the numeral system of  Warlpiri
and Pirah\~{a}  is explained by the fact that the poorness of this
numeral system leads to such results as
 \beq
  \mbox{`many'}+ 1=
\mbox{`many'},   \hspace{5mm}    \mbox{`many'} + 2 = \mbox{`many'},
  \label{Koch_9}
\eeq
 \beq
\mbox{`many'}- 1= \mbox{`many'},   \hspace{5mm}    \mbox{`many'} - 2
= \mbox{`many'},
 \label{Koch_10}
\eeq
 \beq
\mbox{`many'}+   \mbox{`many'} = \mbox{`many'}
 \label{Koch_11}
\eeq
 that are crucial for changing   our outlook on infinity. In
fact, by changing in these relations `many' with $\infty$ we get
relations that are used for working with infinity   in the
traditional calculus:
 \beq
  \infty + 1= \infty,    \hspace{3mm}    \infty + 2 =
\infty, \hspace{3mm} \infty - 1= \infty,    \hspace{3mm}    \infty -
2 = \infty, \hspace{3mm} \infty +
  \infty =   \infty.   \label{Koch_8}
\eeq
 We can see  that numerals `many' and $\infty$ are used in the
same way and we know that in the case of `many' expressions in
(\ref{Koch_9})--(\ref{Koch_11}) are nothing else but the result of
the lack of appropriate numerals for working with finite quantities.
This analogy allows us to conclude that expressions in
(\ref{Koch_8}) used to work  with infinity are also just the result
of the lack of appropriate numerals, in this case  for working with
infinite quantities. As the numeral `many' is not able to represent
the existing richness of finite numbers, the numeral $\infty$ is not
able to represent the richness of the infinite ones.

Notice that it is well known that numeral systems strongly bound the
possibilities to express numbers and to execute mathematical
operations with them. For instance, the Roman numeral system lacks a
numeral   expressing zero. As a consequence, such expressions as V-V
and II-XI in this numeral system are indeterminate forms. The
introduction of the positional numeral system has allowed people to
avoid indeterminate forms of this type and to execute the required
operations easily.

In order to give the possibility to write down more infinite and
infinitesimal numbers, a new numeral system  has been introduced
recently in \cite{Sergeyev,Poland,informatica,UMI}. It allows people
to express a variety of different infinities and infinitesimals, to
perform numerical computations with them,  and to avoid both
expressions of the type (\ref{Koch_8}) and indeterminate forms such
as $\infty-\infty, \frac{\infty}{\infty}, 0\cdot\infty,$ etc.
present in the traditional calculus and related to limits with an
argument tending to $\infty$  or zero.  The   numeral system from
\cite{Sergeyev,Poland,informatica,UMI} has allowed the author to
propose a corresponding computational methodology and to introduce
the Infinity Computer (see the
 patent \cite{Sergeyev_patent}) being a supercomputer  working numerically with a variety
of infinite and infinitesimal numbers. Notice that even though the
new methodology works with infinite and infinitesimal quantities, it
is not related to symbolic computations practiced in  non-standard
analysis (see \cite{Robinson}) and has an applied, computational
character.

    In order to see the place
of the new approach in the historical panorama of ideas dealing with
infinite and infinitesimal, see
\cite{Kauffman,Lolli,Lolli_2,MM_bijection,Sorbi,
Dif_Calculus,first,UMI,Sergeyev_Garro}. In particular, connections
of the new approach with   bijections is studied in
\cite{MM_bijection} and metamathematical investigations on the new
theory and and its non-contradictory  can be found in
\cite{Lolli_2}. The new methodology has been successfully applied
for studying percolation and biological processes (see
\cite{Iudin,Iudin_2,DeBartolo,Biology}), infinite series (see
\cite{Kanovei,Dif_Calculus,Riemann,Zhigljavsky}), hyperbolic
geometry (see \cite{Margenstern,Margenstern_3}), fractals (see
\cite{Iudin,Iudin_2,chaos,Menger,Biology}), numerical
differentiation and optimization (see
\cite{DeLeone,Num_dif,Zilinskas}), the first Hilbert problem, Turing
machines, and lexicographic ordering (see
\cite{first,Sergeyev_Garro,Sergeyev_Garro_2,Sergeyev_Garro_3,medals}),
cellular automata (see \cite{DAlotto,DAlotto_3,DAlotto_2}), ordinary
differential equations (see \cite{ODE,Sergeyev_AIP_2}), etc.

In this paper, the new numeral system  and the respective
computational methodology are used   to study the Koch snowflake.
Both the system and the methodology are based on the introduction in
the process of computations of a new numeral, \G1, called
\textit{grossone}. It is defined as   the infinite integer being the
number of elements of the set, $\mathbb{N}$, of natural
numbers\footnote{Notice that nowadays not only positive integers but
also zero is frequently included in $\mathbb{N}$. However, since
zero has been invented significantly later than positive integers
used for counting objects, zero is not include in $\mathbb{N}$ in
this text.}. Symbols used traditionally to deal with infinite and
infinitesimal quantities (e.g., $\infty$, Cantor's $\omega$,
$\aleph_0, \aleph_1, ...$, etc.) are not used together with \G1.
Similarly, when the positional numeral system and the numeral 0
expressing zero had been introduced, symbols I,   IV, VI,  XIII, and
other symbols from the Roman numeral system had been substituted by
the respective Arabic symbols.

 The
numeral \G1 allows one to express a variety of  numerals
representing different infinities and infinitesimals, to order them,
and to execute numerical computations with all of them in a handy
way. For example, for $\mbox{\ding{172}}$ and
$\mbox{\ding{172}}^{3.1}$ (that are examples of infinities) and
$\mbox{\ding{172}}^{-1}$ and $\mbox{\ding{172}}^{-3.1}$ (that are
examples of infinitesimals) it follows
 \beq
 0 \cdot \mbox{\ding{172}} =
\mbox{\ding{172}} \cdot 0 = 0, \hspace{3mm}
\mbox{\ding{172}}-\mbox{\ding{172}}= 0,\hspace{3mm}
\frac{\mbox{\ding{172}}}{\mbox{\ding{172}}}=1, \hspace{3mm}
\mbox{\ding{172}}^0=1, \hspace{3mm} 1^{\mbox{\tiny{\ding{172}}}}=1,
\hspace{3mm} 0^{\mbox{\tiny{\ding{172}}}}=0,
 \label{3.2.1}
       \eeq
\[
 0 \cdot \mbox{\ding{172}}^{-1} =
\mbox{\ding{172}}^{-1} \cdot 0 = 0, \hspace{5mm}
\mbox{\ding{172}}^{-1} > 0, \hspace{5mm} \mbox{\ding{172}}^{-3.1} >
0, \hspace{5mm} \mbox{\ding{172}}^{-1}-\mbox{\ding{172}}^{-1}= 0,
\]
\[
\frac{\mbox{\ding{172}}^{-1}}{\mbox{\ding{172}}^{-1}}=1,
\hspace{4mm}
\frac{5+\mbox{\ding{172}}^{-3.1}}{\mbox{\ding{172}}^{-3.1}}=5\mbox{\ding{172}}^{3.1}+1,
\hspace{4mm} (\mbox{\ding{172}}^{-1})^0=1, \hspace{4mm}
\mbox{\ding{172}} \cdot \mbox{\ding{172}}^{-1} =1,
       \]
       \[
 \mbox{\ding{172}} \cdot \mbox{\ding{172}}^{-3.1}
=\mbox{\ding{172}}^{-2.1}, \hspace{2mm}
\frac{\mbox{\ding{172}}^{3.1}+4\G1}{\mbox{\ding{172}}}=\G1^{2.1}+4,
\hspace{3mm}
\frac{\mbox{\ding{172}}^{3.1}}{\mbox{\ding{172}}^{-3.1}}=\G1^{6.2},
  \]
       \[
        (\mbox{\ding{172}}^{3.1})^0=1, \hspace{3mm}
\mbox{\ding{172}}^{3.1} \cdot \mbox{\ding{172}}^{-1} =\G1^{2.1},
\hspace{3mm} \mbox{\ding{172}}^{3.1} \cdot \mbox{\ding{172}}^{-3.1}
=1.
       \]

  It follows from (\ref{3.2.1}) that
 a finite number $b$ can be
represented in this numeral system simply as
$b\mbox{\ding{172}}^0=b$, since $\mbox{\ding{172}}^0=1$,  where the
numeral $b$ itself can be written down by any convenient numeral
system used to express finite numbers. The simplest infinitesimal
numbers are represented by numerals having only negative finite
powers of \G1 (e.g., the number
$5.1\mbox{\ding{172}}^{-1.2}\small{+} 6.8\mbox{\ding{172}}^{-20.3}$
consists of two infinitesimal parts, see also examples above).
Notice that all infinitesimals are not equal to zero. For instance,
$\mbox{\ding{172}}^{-3.1} =\frac{1}{\mbox{\ding{172}}^{3.1}}$ is
positive because it is the result of division between two positive
numbers.

\begin{table}[t]
\caption{Measuring infinite sets using   \G1-based numerals allows
one in certain cases to obtain more precise answers in comparison
with the traditional cardinalities, $\aleph_0$ and~$\mathcal{C}$, of
Cantor.}
\begin{center} \scriptsize \label{table1}
\begin{tabular}{@{\extracolsep{\fill}}|c|c|c| }\hline
    &   &      \vspace{-2mm}\\
Description of sets & Cardinality & Number of elements  \\
  &   &      \vspace{-2mm} \\
 \hline
  &   &     \vspace{-2mm}   \\
the set of natural numbers  $\mathbb{N}$   &  countable, $\aleph_0$ & \G1     \\
  &   &     \vspace{-2mm}   \\
  $\mathbb{N} \setminus \{ 3, 5, 10, 23 \} $   &  countable, $\aleph_0$ & \G1-4     \\
  &   &     \vspace{-2mm}   \\
 the set of even numbers $\mathbb{E}$  &  countable, $\aleph_0$ & $\frac{\G1}{2}$    \\
 &   &     \vspace{-2mm}   \\
 the set of odd numbers  $\mathbb{O}$  &  countable, $\aleph_0$ & $\frac{\G1}{2}$    \\
 &   &     \vspace{-2mm}   \\
the set of square natural numbers $\mathbb{G} = \{ x : x= n^2, n \in \mathbb{N}, x \in \mathbb{N} \}$   &  countable, $\aleph_0$ & $ \lfloor \sqrt{\G1} \rfloor$      \\
 &   &     \vspace{-2mm}   \\
 the set of integer numbers $\mathbb{Z}$    &  countable, $\aleph_0$ & 2\G1+1   \\
 &   &     \vspace{-2mm}\\
 the set of pairs of natural numbers $\mathbb{P}  = \{ (p,q) : p   \in
\mathbb{N}, q \in \mathbb{N} \}$   &  countable, $\aleph_0$ &  $\G1^2$   \\
  &   &     \vspace{-2mm}\\
 the set of numerals  $\mathbb{Q}'  = \{     -\frac{p}{q}, \,\, \frac{p}{q} : p   \in
\mathbb{N}, q \in \mathbb{N} \} $ &  countable, $\aleph_0$ &  $2\G1^2$   \\
  &   &     \vspace{-2mm}\\
 the set of numerals  $\mathbb{Q}  = \{ 0,   -\frac{p}{q}, \,\, \frac{p}{q} : p   \in
\mathbb{N}, q \in \mathbb{N} \} $ &  countable, $\aleph_0$ &  $2\G1^2+1$   \\
  &   &     \vspace{-2mm}\\
 the set of numerals $A_2$  &  continuum, $\mathcal{C}$ & $2^{\mbox{\tiny{\ding{172}}}}$     \\
 &   &     \vspace{-2mm}   \\
 the set of numerals   $A'_2  $  &  continuum, $\mathcal{C}$ & $2^{\mbox{\tiny{\ding{172}}}}+1$     \\
 &   &     \vspace{-2mm}   \\
the set of numerals $A_{10}$  &  continuum, $\mathcal{C}$ & $10^{\mbox{\tiny{\ding{172}}}}$    \\
 &   &     \vspace{-2mm}   \\
 the set of numerals $C_{10}$  &  continuum, $\mathcal{C}$ & $2 \cdot 10^{\mbox{\tiny{\ding{172}}}}$    \\
\hline
\end{tabular}
\end{center}
\label{Koch_table}
\end{table}

In the context of the present paper it is important that in
comparison to the traditional mathematical tools used to work with
infinity the new numeral system allows one to obtain more precise
answers in certain cases.   For instance, Tab.~\ref{Koch_table}
compares results obtained by the traditional Cantor's cardinals and
the new numeral system with respect to the measure of a great dozen
of infinite sets (for a detailed discussion regarding the results
presented in Tab.~\ref{Koch_table} and for more examples dealing
with infinite sets see
\cite{Lolli_2,MM_bijection,first,Lagrange,Sergeyev_Garro}). Notice,
that in $\mathbb{Q}$ and $\mathbb{Q}'$ we calculate different
numerals and not numbers. For instance, numerals $\frac{3}{1}$ and
$\frac{6}{2}$ have been counted two times even though they represent
the same number~3. Then, four sets of numerals having the
cardinality of continuum are  shown in Tab.~\ref{Koch_table}. Among
them we denote by $A_2$  the set of numbers   $x \in [0,1)$
expressed in the binary positional numeral system,  by  $A'_2 $ the
set being the same as  $A_2  $ but with $x$ belonging to the closed
interval $[0,1]$, by $A_{10}$ the set of numbers $x \in [0,1)$
expressed in the decimal positional numeral system, and finally we
have the set $C_{10} = A_{10} \cup B_{10}$, where $B_{10}$ is the
set of numbers $x \in [1,2)$ expressed in the decimal positional
numeral system. It is worthwhile to notice also that  \G1-based
numbers present in Tab.~\ref{Koch_table} can be ordered as follows
\[
\lfloor \sqrt{\G1} \rfloor < \frac{\G1}{2} < \G1-4  < \G1 < 2\G1 <
2\G1+1 <
 \]
 \[
  \G1^2 < 2\G1^2+1 < 2^{\mbox{\tiny{\ding{172}}}} <
2^{\mbox{\tiny{\ding{172}}}}+1 < 10^{\mbox{\tiny{\ding{172}}}} < 2
\cdot 10^{\mbox{\tiny{\ding{172}}}}.
\]

It can be seen from Tab.~\ref{Koch_table} that Cantor's
cardinalities say only whether a set is countable or uncountable
while the \G1-based numerals allow us to express the exact number of
elements of the infinite sets. However, both numeral systems -- the
new one and the numeral system of infinite cardinals -- do not
contradict one another. Both   numeral systems provide correct
answers, but their answers have  different accuracies. By using an
analogy from physics we can say that the lens of our new `telescope'
used to observe infinite sets is stronger and where Cantor's
`telescope' allows one to distinguish just two dots (countable sets
and the continuum) we are able to see many different dots (infinite
sets having different number of elements).

\section{Quantitative characteristics of the Koch snowflake at infinity}

As was mentioned above, \G1 can be successfully used for various
purposes related to studying infinite and infinitesimal objects, in
particular,  for indicating positions of elements in infinite
sequences (see, e.g., \cite{Sergeyev_Garro,Sergeyev_Garro_2,UMI})
and for working with divergent series (see
\cite{Kanovei,Dif_Calculus,Riemann,Zhigljavsky}). Both topics will
help us to study the Koch snowflake at infinity. Let us first
compare how infinite sequences are defined from the traditional
point of view and from the new one.

The traditional definition   is very simple:  An infinite sequence
$\{b_n\},$  where for all $n \in \mathbb{N}$ elements $b_n$ belong
to a set $B$ is defined as a function having as its domain the set
$\mathbb{N}$ and as its codomain the set~$B$. Let us see now how
this definition can be reformulated using the new methodology and
\G1-based numerals. Remind that \G1 has been introduced as the
number of elements of the set   $\mathbb{N}$. Thus, due to the
traditional definition given above, any sequence having $\mathbb{N}$
as the domain  has \ding{172} elements.

In its turn, the notion of a subsequence is introduced traditionally
as a sequence from which some of its elements have been deleted. In
cases where both the original sequence and the obtained subsequence
are infinite, in spite of the fact that some elements were excluded,
the traditional fashion does not allow us to record in some way that
the obtained infinite subsequence has less elements than the
original infinite sequence. In the new fashion there is such a
possibility. Having a sequence with \G1 elements exclusion of $k$
elements from it gives a subsequence having $\G1-k < \G1$  elements.
For instance, in (\ref{Koch_2}) the first infinite sequence has \G1
elements and the second one \G1-2 elements:
 \beq
\underbrace{1,2,3,4,\hspace{1mm}  \ldots \hspace{1mm}
\mbox{\ding{172}}-2,\hspace{1mm}
 \mbox{\ding{172}}-1,
\mbox{\ding{172}}}_{\mbox{\ding{172} elements}}, \hspace{5mm}
\underbrace{ 4,5,6, \hspace{1mm}  \ldots \hspace{1mm}
\mbox{\ding{172}}-2,\hspace{1mm}
 \mbox{\ding{172}}-1, \G1, \G1+1}_{\mbox{\ding{172}-2 elements}}.
  \label{Koch_2}
 \eeq
Thus, the numeral system using \ding{172} allows  us to observe not
only the starting but also the \textit{ending} elements of infinite
processes, if the respective elements are expressible in this
numeral system. This fact is   important in connection with fractals
because it allows us to distinguish different fractal objects after
an infinite number of steps of their construction. Another useful
observation consists of the fact that, since the number of elements
of any   sequence (finite or infinite) is less or equal
to~\ding{172},  any sequential process can have at maximum
\ding{172} steps  (see  \cite{Sergeyev_Garro}).

Let us see now what the new approach allows us to say with respect
to divergent series. In particular, the situations where it is
necessary to sum up an infinite number of infinitesimal numbers will
be  of our primary interest (for a detailed discussion see
\cite{Kanovei,Dif_Calculus,Riemann,Zhigljavsky}). This issue is
important for us since in (\ref{Koch_1}) the perimeter $P_n$ of the
snowflake studied traditionally goes to infinity. Grossone-based
numerals allow us to express not only different finite numbers but
also different infinite numbers so, such expressions as
$S_1=b_1+b_2+\ldots$ or $S_1=\sum_{i=1}^{\infty}b_i$ become
unprecise since the number of addends in $S_1$ is not indicated
explicitly. If we use again the analogy with  Warlpiri and
Pirah\~{a} then we can say that the record $\sum_{i=1}^{\infty}b_i$
can be interpreted as $\sum_{i=1}^{many}b_i$. Notice that in the
finite sums the situation is the same: it is not sufficient to say
that the number of summands is finite, it is necessary to define
explicitly their number.

The new approach gives the possibility to add finite, infinite, and
infinitesimal values in a handy way, the number of summands can be
finite or infinite, and results of addition can be finite, infinite,
and infinitesimal in dependence on the sort and number of addends.
To illustrate this assertion let us consider a few examples. First,
it becomes possible to compute the sum of all natural numbers from 1
to \ding{172}  as follows
 \beq
1+2+3+ \ldots + (\mbox{\ding{172}}-1) + \mbox{\ding{172}} =
\sum_{i=1}^{\mbox{\tiny{\ding{172}}}} i =
\frac{\mbox{\ding{172}}}{2}(1 + \mbox{\ding{172}})=
0.5\mbox{\ding{172}}^{2}\mbox{\small{+}}0.5\mbox{\ding{172}}.
 \label{Koch_6}
 \eeq
 The following sum of
infinitesimals where each summand is \ding{172} times less than the
corresponding item of (\ref{Koch_6}) can be also computed easily
 \[
\mbox{\ding{172}}^{-1}+2\mbox{\ding{172}}^{-1}+ \ldots +
(\mbox{\ding{172}}-1)\cdot\mbox{\ding{172}}^{-1} +
\mbox{\ding{172}}\cdot\mbox{\ding{172}}^{-1} =
 \]
\beq \sum_{i=1}^{\mbox{\tiny{\ding{172}}}} i\mbox{\ding{172}}^{-1} =
\frac{\mbox{\ding{172}}}{2}(\mbox{\ding{172}}^{-1} + 1)=
0.5\mbox{\ding{172}}^{1}\mbox{\small{+}}0.5.
 \label{Koch_7}
 \eeq
As expected, the obtained number, $0.5\mbox{\ding{172}}^{1}+0.5$ is
\ding{172} times less than the result obtained in (\ref{Koch_6}).
Notice that this example shows, in particular, that sum of
infinitely many infinitesimals can be infinite.

Then, in the same way as it happens in situations where the number
of summands is finite, the following examples show that   smaller or
larger number of summands changes the result (cf. (\ref{Koch_6}),
(\ref{Koch_7}))
 \[
\sum_{i=1}^{\mbox{\tiny{\ding{172}}}-1} i =
\frac{\mbox{\ding{172}}-1}{2}(1 + \mbox{\ding{172}}-1)=
0.5\mbox{\ding{172}}^{2}\mbox{\small{-}}0.5\mbox{\ding{172}},
 \]
\[
\sum_{i=1}^{3\mbox{\tiny{\ding{172}}}^2} i\mbox{\ding{172}}^{-1} =
 \frac{3\mbox{\ding{172}}^2}{2}(\mbox{\ding{172}}^{-1} +
3\G1)= 4.5\mbox{\ding{172}}^{3}+1.5\G1^2.
\]
Notice that sums can have more than \G1 addends if it is not
required to execute the operation of addition by a successive adding
summands, i.e., the summation can be   done in parallel. However, if
in a particular application there exists a restriction that the
required summation should be executed sequentially, then, since any
sequential process  cannot have more than \G1 steps, the sequential
process of the summation cannot have more than \G1 addends.

We are ready now to return to the Koch snowflake and to study it at
infinity using \G1-based numerals. It can be seen from
Fig.~\ref{Koch_fig} that at each iteration each side of the
snowflake is substituted by 4 new sides having the length of one
third of the segment that has been substituted. Thus, if we indicate
as $N_n, n \ge 1,$ the number of segments of the snowflake and    as
$L_n$ their length at the $n$-th iteration then
 \beq N_n = 4N_{n-1} = 3 \cdot 4^n, \hspace{1cm}
 L_n=\frac{1}{3}L_{n-1}= \frac{l}{3^n}, \hspace{1cm} n > 1,
 \label{Koch_3}
 \eeq
where $l$ is the length of each side of the original triangle from
Fig.~\ref{Koch_fig}. As a
 result, the perimeter $P_n$   of the Koch snowflake is calculated
as follows
 \beq
  P_n = N_n \cdot L_n = \frac{4}{3} N_{n-1} \cdot L_{n-1} = \frac{4}{3}
  P_{n-1} =
\frac{4^n}{3^{n-1}}l. \label{Koch_4}
 \eeq
Therefore, if  we start our computations from the original triangle
from Fig.~\ref{Koch_fig}, after \G1 steps we have the snowflake
having the infinite number of segments $N_{\tiny{\G1}} =   3 \cdot
4^{\tiny{\G1}}$. Each of the segments has the infinitesimal length
$L_{\tiny{\G1}} = \frac{1}{3^{\tiny{\G1}}}l$. In order to calculate
the perimeter, $P_{\tiny{\G1}}$, of the snowflake we should multiple
the infinite number $N_{\tiny{\G1}}$ and the infinitesimal number
$L_{\tiny{\G1}}$. Thus the perimeter is
 \[
 P_{\tiny{\G1}} = N_{\tiny{\G1}}   \cdot  L_{\tiny{\G1}} = 3 \cdot
   4^{\tiny{\G1}} \cdot
\frac{1}{3^{\tiny{\G1}}}l = \frac{4^{\tiny{\G1}}}{3^{\tiny{\G1}-1}}l
 \]
and it is infinite. Analogously, in case we have executed \G1-1
steps we have the infinite perimeter $P_{\tiny{\G1-1}} =
\frac{4^{\tiny{\G1}-1}}{3^{\tiny{\G1}-2}}l$. Since the new numeral
systems allows us to execute easily arithmetical operations with
infinite numbers, we can   divide the obtained two infinite numbers,
$P_{\tiny{\G1}}$ and $P_{\tiny{\G1}-1}$, one by another and to
obtain as the result the finite number that is in a complete
agreement with (\ref{Koch_4})
 \[
 \frac{P_{\tiny{\G1}}}{P_{\tiny{\G1-1}}}  = \frac{
\frac{4^{\tiny{\G1}}}{3^{\tiny{\G1}-1}}l}{\frac{4^{\tiny{\G1}-1}}{3^{\tiny{\G1}-2}}l}
= \frac{4}{3}.
 \]
 The difference of the two perimeters can also be calculated easily
 \[
 P_{\tiny{\G1}}-P_{\tiny{\G1-1}}   =
\frac{4^{\tiny{\G1}}}{3^{\tiny{\G1}-1}}l-\frac{4^{\tiny{\G1}-1}}{3^{\tiny{\G1}-2}}l
= \frac{4^{\tiny{\G1}-1}}{3^{\tiny{\G1}-2}}l \left(  \frac{4}{3}-1
\right) = \frac{4^{\tiny{\G1}-1}}{3^{\tiny{\G1}-1}}l
 \]
and it results to be infinite.

In case the infinite number of steps $n=0.5\G1$, it follows that the
infinite perimeter is $P_{0.5\tiny{\G1}} =
\frac{4^{0.5\tiny{\G1}}}{3^{0.5\tiny{\G1}-1}}l$ and the operation of
division of two infinite numbers gives us as the result also an
infinite number that can be calculated precisely:
 \[
 \frac{P_{\tiny{\G1}}}{P_{0.5\tiny{\G1}}}  = \frac{
\frac{4^{\tiny{\G1}}}{3^{\tiny{\G1}-1}}l}{\frac{4^{0.5\tiny{\G1}}}{3^{0.5\tiny{\G1}-1}}l}
= \left(\frac{4}{3}\right)^{0.5\tiny{\G1}}.
 \]
Thus we can distinguish now in a precise manner that the infinite
perimeter $P_{\tiny{\G1}}$ is infinitely times  longer than the
infinite perimeter $P_{\tiny{0.5\G1}}$.

We can distinguish also at infinity  the snowflakes having different
initial generators. As we have already seen, starting from the
original triangle after \G1 steps the snowflake has $   3 \cdot
4^{\tiny{\G1}}$ segments and each of them has the infinitesimal
length $L_{\tiny{\G1}} = \frac{1}{3^{\tiny{\G1}}}l$ and  the
  perimeter of the snowflake  is $P_{\tiny{\G1}} =
\frac{4^{\tiny{\G1}}}{3^{\tiny{\G1}-1}}l$. If we start from the
initial configuration shown at Fig.~\ref{Koch_fig}, c) and also
execute \G1 steps then the resulting snowflake will have
$N_{\tiny{\G1}+2} =   3 \cdot 4^{\tiny{\G1}+2}$ segments, each of
them will have the infinitesimal length $L_{\tiny{\G1}+2} =
\frac{1}{3^{\tiny{\G1}+2}}l$ and the perimeter of the snowflake will
be $P_{\tiny{\G1}+2} = \frac{4^{\tiny{\G1}+2}}{3^{\tiny{\G1}+1}}l$.
Thus, this snowflake will have infinitely more segments than the one
started from the original triangle. More precisely, this infinite
difference is equal to
\[
N_{\tiny{\G1}+2} - N_{\tiny{\G1}}=   3 \cdot 4^{\tiny{\G1}+2}  - 3
\cdot 4^{\tiny{\G1}} = 3 \cdot 4^{\tiny{\G1}} (4^2-1) = 45 \cdot
4^{\tiny{\G1}}.
\]
The lengths of the segments in both snowflakes are infinitesimal
and, in spite of the fact that their difference is also
infinitesimal, it   can   be calculated precisely as follows
\[
L_{\tiny{\G1}} - L_{\tiny{\G1}+2}= \frac{1}{3^{\tiny{\G1}}}l -
\frac{1}{3^{\tiny{\G1}+2}}l =
\frac{1}{3^{\tiny{\G1}}}l\left(1-\frac{1}{3^2} \right) =
\frac{8l}{3^{\tiny{\G1}+2}}.
\]

Let us see now what happens with the area of the snowflake at
infinity. As it can be seen from Fig.~\ref{Koch_fig}, at each
iteration $n$ a new triangle is added at each side of the snowflake
built at iteration $n-1$ and, therefore, the number of new
triangles, $T_n$, is equal to
 \[
 T_n = N_{n-1} =  3 \cdot 4^{n-1}.
 \]
The area, $a_n$, of each triangle added at $n$-th iteration is
$\frac{1}{9}$ of each triangle added during the iteration $n-1$
 \[
 a_n = \frac{a_{n-1}}{9} =  \frac{a_0}{9^{n}},
 \]
where $a_0$ is the area of the original triangle (see
(\ref{Koch_1})). Therefore, the whole new area added to the
snowflake is
  \beq
 T_n   a_n =    3 \cdot 4^{n-1} \cdot  \frac{a_0}{9^{n}} =  \frac{3}{4}
 \cdot\left(\frac{4}{9}\right)^n \cdot
 a_0 = \frac{ a_0}{3}
 \cdot\left(\frac{4}{9}\right)^{n-1}
 \label{Koch_5}
 \eeq
and the complete area, $A_n$, of the snowflake at the $n$-th
iteration is
 \[
  A_n = a_0 + \sum_{i=1}^{n} T_n   a_n =
  a_0 \left(1 +  \frac{1}{3}\sum_{i=1}^{n}\left(\frac{4}{9}\right)^{i-1}
  \right) = a_0 \left(1 +  \frac{1}{3}\sum_{i=0}^{n-1}\left(\frac{4}{9}\right)^{i}
  \right) =
  \]
  \[
a_0 \left(1 +  \frac{1}{3} \cdot
\frac{1-\left(\frac{4}{9}\right)^{n}}{1-\frac{4}{9}} \right)  = a_0
\left(1 +  \frac{3}{5} \left(1-\left(\frac{4}{9}\right)^{n}
  \right)\right) = \frac{a_0 }{5} \left(8-3\left(\frac{4}{9}\right)^{n}
  \right).
  \]
Traditionally, the limit of the area is considered and the result
 \[
  \lim_{n \rightarrow
\infty} A_n = \frac{8}{5}a_0
 \]
is obtained. Thanks to \G1-based numerals we are able now to work
with infinitesimals easily (see
\cite{Dif_Calculus,Riemann,Zhigljavsky} for more results on summing
up infinitesimals, divergent series, etc.) and to observe the
infinitesimal difference of the areas of the snowflakes at infinity.
For example, after \G1-1 and \G1  iterations the difference between
the areas $ A_{\tiny{\G1}-1}$ and $ A_{\tiny{\G1}}$ is
 \[
A_{\tiny{\G1}} -   A_{\tiny{\G1}-1}  =    \frac{a_0 }{5}
\left(8-3\left(\frac{4}{9}\right)^{{\tiny{\G1}}} \right)   -
\frac{a_0 }{5}
  \left(8-3\left(\frac{4}{9}\right)^{{\tiny{\G1}}-1} \right) =
 \frac{a_0}{3}  \left(\frac{4}{9}\right)^{{\tiny{\G1-1}}} > 0.
 \]
This number is infinitesimal and it   perfectly corresponds to the
general formula~(\ref{Koch_5}).

 \section{Conclusion}

Traditional approaches  for   studying dynamics of fractal processes
very often are not able to give their quantitative characteristics
at infinity and, as a consequence, use limits to overcome this
difficulty. The Koch snowflake in fact is defined as the limit
object and, in the limit, its perimeter goes to infinity and, at the
same time, its limit area is finite. As a consequence,
 questions regarding   quantitative
characteristics of  the snowflake at infinity very often remain
without any answer if traditional mathematical tools are used.
Moreover, we cannot distinguish at infinity snowflakes starting from
different initiators even though they are different for any fixed
finite value of the generation step $n$.

In this paper, it has been shown that recently introduced \G1-based
numerals  give the possibility to work with different infinities and
infinitesimals numerically and to establish the presence of
\textit{infinitely many different} snowflakes at infinity instead of
a unique snowflake observed traditionally. For different infinite
values of generation steps $n$ it becomes possible to obtain their
exact quantitative characteristics instead of traditionally made
qualitative declarations saying that limits under consideration go
to infinity.

 In particular, it becomes   possible
to compute the exact infinite number, $N_n$, of sides of the
snowflake, the exact infinitesimal length, $L_n$, of each side for a
given infinite $n$
 and to calculate the exact infinite  value of the
perimeter, $P_n$, of the Koch snowflake. It has been shown that the
sum of  infinitely many   infinitesimal sides gives as the result
the infinite perimeter $P_n$  of the Koch snowflake.  As a
consequence, it has been also shown that for infinite $ k > n$ it
follows that infinitesimal $L_k < L_n$, infinite  $N_k  > N_n$, $P_k
> P_n$ and the exact difference $P_k - P_n$ can be calculated. The
areas $A_n$ and $A_k$ can be also computed exactly for infinite
values $ k$ and $n$. If   $ k
> n$ then it follows that  $A_k > A_n$, the  difference $A_k - A_n$ is infinitesimal
and it can also be calculated exactly.

Moreover, it has been shown the importance of the initial conditions
in the processes of the construction of the Koch snowflake. If we
consider one process of the construction of the snowflake starting,
e.g., from the original triangle and   the initiator of the second
process is a result of first $k$ steps from the original triangle
then after the same infinite number of steps, $n$, the two resulting
snowflakes will be different and it is possible to calculate their
exact perimeters, areas, etc.

\bibliographystyle{plain}
\bibliography{XBib_Koch}
\end{document}